\input amssym.tex
\magnification 1200
\hsize = 14.5cm
\hoffset -0.5cm

\font\Bbb=msbm10
\def\BBB#1{\hbox{\Bbb#1}}
\font\Frak=eufm10
\def\frak#1{{\hbox{\Frak#1}}}
\font\smallFrak=eufm10 scaled 800
\def\smallfrak#1{{\hbox{\smallFrak#1}}}

\hyphenation{Bor-cherds}
\hyphenation{pre-print}


\def\Ldg{1.1}
\def\Ldk{1.2}
\def\Ldd{1.3}

\def\pairDK{{1.4}}

\def\cvp{{3.1}}
\def\hvp{{3.2}}

\def\gamo{{4.1}}
\def\ko{{4.2}}
\def\ka{{4.3}}
\def\gr{{4.4}}
\def\da{{4.5}}
\def\tdd{{4.6}}

\def\drab{{5.1}}
\def\wdoa{{5.2}}
\def\dpz{{5.3}}
\def\Ydab{{5.4}}
\def\dwd{{5.5}}
\def\expr{{5.6}}
\def\barhc{{5.7}}
\def\embd{{5.8}}
\def\relV{{5.9}}
\def\relM{{5.10}}
\def\relVY{{5.11}}
\def\cab{{5.12}}
\def\cenc{{5.13}}
\def\doo{{5.14}}
\def\wdaa{{5.15}}
\def\cencp{{5.16}}

\def\FA{{3.1}}
\def\emb{{5.1}}
\def\Memb{{5.2}}
\def\presr{{5.3}}
\def\maind{{5.4}}
\def\irrd{{5.5}}
\def\TDK{{5.6}}


\def\D{{\cal D}}

\def\HVir{{{\cal H}{\cal V}{\mit ir}}}
\def\DK{{{\cal D}{\cal K}}}

\def\d{\partial}
\def\g{{\frak g}}
\def\f{{\frak f}}
\def\sg{{\smallfrak g}}
\def\sf{{\smallfrak f}}

\def\glN{{\it gl}_N}
\def\gln{{\it gl}_N}

\def\slN{{{\it sl}_N}}
\def\sln{{{sl}_N}}
\def\wsl{{\widehat {sl}_N}}
\def\dg{{\dot \g}}
\def\sdg{{\dot \sg}}

\def\wdg{{\widehat \dg}}
\def\swdg{{\widehat \sdg}}
\def\td{\tilde d}
\def\wda{{{\widehat d}_a}}
\def\wL{{\overline L}}
\def\wcL{{{\overline c}_\Vir}}
\def\wCL{{{\overline C}_\Vir}}
\def\df{{\dot \f}}

\def\tf{{\tilde \f}}
\def\barf{{\bar \f}}
\def\sbarf{{\bar \sf}}

\def\om{\omega}
\def\ot{\otimes}
\def\eru{e^{ru}{}}

\def\div{{\rm div}}
\def\gdiv{{\g_\div}}
\def\Ddiv{\D_\div}

\def\sgdiv{{\sg_\div}}

\def\Hyp{{\mit Hyp}}
\def\hyp{{\mit Hyp}}

\def\VH{{V_{\Hyp}^+}}

\def\Hei{{{\cal H}{\mit ei}}}
\def\HV{{\cal VH}}
\def\R{{\cal R}}
\def\K{{\cal K}}

\def\C{\BBB C}
\def\o{{\bf 1}}
\def\Z{\BBB Z}

\def\Id{{\rm Id}}

\def\Der{{{\rm Der}\hbox{\hskip 0.08cm}}}

\def\Vir{{{\cal V}{\mit ir}}}
\def\wVir{{\overline \Vir}}

\def\char{{{\rm char}\hbox{\hskip 0.08cm}}}
\def\tr{{\rm tr}}

\def\arsubset{{\hbox{\hskip 2pt}\subset{\kern-8pt{\lower2.7pt\hbox{$\rightarrow$}}}\hbox{\hskip 3pt}}}

\

\

\centerline
{\bf Representations of toroidal extended affine Lie algebras.}

\centerline{
{\bf Yuly Billig}
\footnote{*}{Research supported by the  Natural Sciences and
Engineering Research Council of Canada.}
}

\

\centerline{School of Mathematics \& Statistics}
\centerline{Carleton University}
\centerline{1125 Colonel By Drive}
\centerline{Ottawa, Ontario, K1S 5B6, Canada}
\centerline{e-mail: billig@math.carleton.ca}

\

\

\

{\bf Abstract.} We show that the representation theory for the toroidal extended affine
Lie algebra is controlled by a VOA which is a tensor product of four VOAs: a sub-VOA $\VH$ of
a hyperbolic lattice VOA, affine $\wdg$ and $\wsl$ VOAs and a Virasoro VOA.
A tensor product of irreducible modules for these VOAs admits the structure of 
an irreducible module for the toroidal extended affine Lie algebra.
We also show that for $N=12$, $\VH$ becomes an exceptional irreducible module for
the extended affine Lie algebra of rank 0.

\


\

\

{\bf 0. Introduction.}

\

In this paper we study representations of toroidal extended affine Lie algebras
using recently developed representation theory for the full toroidal Lie algebras [B3].
Extended affine Lie algebras (EALAs) have been extensively studied during the last decade (see [N], [ABFP], [AABGP]
and references therein). The main features of an extended affine Lie algebra
is that it is graded by a finite root system and possesses a non-degenerate
symmetric invariant bilinear form. 

The construction of toroidal Lie algebras parallels one for affine algebras.
We start with the Lie algebra of maps from an $N+1$-dimensional torus into a finite-dimensional
simple Lie algebra $\dg$. This multi-loop algebra may be written as a tensor product $\R \otimes \dg$
of the algebra $\R$ of Laurent polynomials in $N+1$ variables with $\dg$.
Next, we take the universal central extension $\R \otimes \dg \oplus \K$ of this 
multi-loop algebra and add a Lie algebra of vector fields on the torus, possibly twisted with
a 2-cocycle (see section 1 for details). If we add all vector fields, we get the full toroidal
Lie algebra. However, the full toroidal Lie algebra does not
possess a non-degenerate invariant form, whereas its subalgebra with the divergence zero vector fields
$$ \gdiv = (\R \ot \dg) \oplus \K \oplus \D_\div$$
does. We call this last algebra the toroidal extended affine Lie algebra.

The representation theory of toroidal Lie algebras is best described in the framework
of vertex operator algebras (VOAs).
We prove in this paper that the vertex operator algebra that controls the representation theory of $\gdiv$ is
a tensor product of an affine $\wdg$ VOA, a sub-VOA of a hyperbolic
lattice VOA $V_\hyp^+$, affine $\wsl$ VOA and a Virasoro VOA. By taking a tensor product of irreducible 
modules for each of these VOAs we get irreducible modules for the toroidal extended affine algebra $\gdiv$
(Theorem \irrd). This gives an explicit realization for these modules and allows us to compute their
characters.

The results of this paper can be used to construct irreducible modules for virtually all other EALAs [BL]. 

If we set $\dg = (0)$, we will get representations for the Lie algebra $\D_\div \oplus \K$. 
This Lie algebra still possesses a non-degenerate symmetric invariant bilinear form and may be viewed as
an EALA of rank 0. This Lie algebra plays an important role in magnetic hydrodynamics [VD], [B2].
We establish one curious fact about the representation theory of $\D_\div \oplus \K$.
When $N = 12$, this Lie algebra has an exceptional module with a particularly simple structure.
Only for this value of $N$ a hyperbolic lattice sub-VOA $V_\hyp^+$ can be given a structure of an irreducible
module for $\D_\div \oplus \K$. The character
of this module is given by the $-24$-th power of the Dedekind $\eta$-function
and has nice modular properties.

 This paper is a sequel to [B3], and we will be using constructions and notations from 
that paper. A weaker form of the results of this paper was presented in [B1] (unpublished).

 This paper is structured as follows. In section 1 we review the construction of toroidal Lie 
algebras. In sections 2 and 3 we discuss vertex operator algebras that we will need in order
to build the representation theory of the toroidal EALAs. In section 4 we recall the results
on the representations of the full toroidal Lie algebras. The final section 5 contains the main 
results of this paper, there we develop the representation theory for toroidal extended affine
Lie algebras.

\

 {\bf Acknowledgments:} This work has been completed during my stay at the University of Sydney.
I thank the School of Mathematics and Statistics, University of Sydney, for the warm hospitality.

\

\

{\bf 1. Toroidal Lie algebras.}

\

Toroidal Lie algebras are the natural multi-variable generalizations of
affine Lie algebras. 
In this review of the toroidal Lie algebras we follow the work [BB].
Let $\dg$ be a simple finite-dimensional Lie algebra
over $\C$ with a non-degenerate invariant bilinear form $( \cdot | \cdot )$
and let $N \geq 1$ be an integer. 
We consider the Lie algebra $\R \ot \dg$
of maps of an $N+1$ dimensional torus into $\dg$, where  
$\R = \C [t_0^\pm, t_1^\pm, \ldots, t_N^\pm]$ is the algebra of functions 
on a torus (in the Fourier basis). The universal central extension
of this Lie algebra may be described by means of the following construction which is
due to Kassel [Kas]. Let $\Omega_\R$ be the space of $1$-forms on a torus:
$\Omega_\R = \mathop\oplus\limits_{p = 0}^N \R dt_p$. We will choose the
forms $\{ k_p = t_p^{-1} dt_p \, | \, p = 0, \ldots, N \}$ as a basis of this
free $\R$ module. There is a natural map $d$ from the space of functions
$\R$ into $\Omega_\R$: $d(f) = \sum\limits_{p = 0}^N {\d f \over \d t_p} dt_p
= \sum\limits_{p = 0}^N t_p {\d f \over \d t_p} k_p$. The center $\K$ for the
universal central extension $(\R \ot \dg) \oplus \K$ of $\R \ot \dg$
is realized as
$$ \K = \Omega_\R / d(\R) , $$
and the Lie bracket is given by the formula
$$[f_1(t) g_1, f_2(t) g_2] = f_1(t) f_2(t) [g_1, g_2] + (g_1| g_2) f_2 d(f_1).$$
Here and in the rest of the paper we will denote elements of $\K$ by the
same symbols as elements of $\Omega_\R$, keeping in mind the canonical
projection $\Omega_\R \rightarrow \Omega_\R / d(\R)$.

Just as in affine case, we add to $(\R \ot \dg) \oplus \K$ the algebra $\D$ of outer derivations 
$$\D = \mathop\oplus\limits_{p=0}^N \R d_p,$$
where $ d_p = t_p {\d \over \d t_p}$.
We will use the multi-index notations,  
$r = (r_0, r_1, \ldots, r_N)$, 
$t^r = t_0^{r_0} t_1^{r_1} \ldots t_N^{r_N}$.

The natural action of $\D$ on $\R \ot \dg$ 
$$[t^r d_a, t^m g] = m_a t^{r+m} g \eqno{(\Ldg)}$$
uniquely extends to the action on the universal central extension
$(\R \ot \dg) \oplus \K$ by
$$[t^r d_a, t^m k_b] = m_a t^{r+m} k_b + \delta_{ab}
\sum\limits_{p=0}^N r_p t^{r+m} k_p . \eqno{(\Ldk)}$$
This corresponds to the Lie derivative action of the vector fields on 1-forms.

It turns out that there is still an extra degree of freedom in defining
the Lie algebra structure on $(\R \ot \dg) \oplus \K \oplus \D$.
The Lie bracket on $\D$ may be twisted with a $\K$-valued 2-cocycle:
$$[t^r d_a, t^m d_b] = m_a t^{r+m} d_b - r_b t^{r+m} d_a
+ \mu m_a r_b \sum\limits_{p=0}^N m_p t^{r+m} k_p . \eqno{(\Ldd)}$$
As a result we get a family of full toroidal Lie algebras 
$$\g(\mu) = (\R \ot \dg) \oplus \K \oplus \D.$$

Note that after adding the algebra of derivations $\D$, the center
of the toroidal Lie $\g$ becomes finite-dimensional with the basis 
$\{ k_0, k_1, \ldots, k_N \}$. This can be seen from the action
(\Ldk) of $\D$ on $\K$,  which is non-trivial.

The toroidal Lie algebra $\g(\mu) = (\R \otimes \dg) \oplus \K \oplus \D$
has an important subalgebra $\gdiv (\mu)$ that has divergence free vector fields
as the derivation part:
$$ \gdiv = \gdiv(\mu) = (\R \otimes \dg) \oplus \K \oplus \Ddiv,$$
where
$$\Ddiv = \left\{ \sum_{p=0}^N f_p(t) d_p \quad \bigg| \quad \sum_{p=0}^N t_p
{\d f_p \over \d t_p} = 0 \right\} .$$
The expression $i \sum\limits_{p=0}^N t_p  {\d f_p \over \d t_p}$ becomes
the divergence of a vector field in the angular coordinates $(x_0, \ldots, x_N)$
on a torus, where $t_j = e^{i x_j}$.


The importance of this subalgebra is explained by the fact that unlike the full
toroidal Lie algebra, $\gdiv(\mu)$ is an extended affine Lie algebra [BGK],
i.e., $\gdiv (\mu)$ has a non-degenerate symmetric invariant bilinear form. The
restrictions of this form to both $\R \otimes \dg$ and to $\Ddiv \oplus \K$ are
non-degenerate:
$$(t^r g_1 | t^m g_2 ) = \delta_{r,-m} (g_1 | g_2) , \quad g_1, g_2 \in \dg,$$
while the vector fields pair with the $1$-forms:
$$\big( \sum_{p=0}^N a_p t^r d_p | t^m k_q \big) = \delta_{r,-m} a_q . 
\eqno{(\pairDK)}$$
One can see that the above formula is ill-defined for the full $\D$,
since $d(t^m) = \sum\limits_{q=0}^N m_q t^m k_q$, being zero in $\K$, 
must be in the kernel
of the form. For the subalgebra $\Ddiv$ this is precisely the case since
$$ \big( \sum_{p=0}^N a_p t^r d_p |  \sum_{q=0}^N r_q t^{-r} k_q \big) = 
\sum\limits_{q=0}^N a_q r_q = 0.$$

All other values of the bilinear form are trivial:
$$(\R \otimes \dg | \Ddiv \oplus \K) = 0, \quad
(\Ddiv | \Ddiv) = 0, \quad (\K | \K) = 0 .$$
It is easy to verify that the resulting symmetric bilinear form is
invariant and non-degenerate.

 The 2-cocycle that we consider here was originally introduced in [EM] and
was denoted $\mu \tau_1$ in [BB], [B3].
There is another cocycle $\tau_2$ used in [L], [BB], [B3]. However the cocycle
$\tau_2$ vanishes on $\D_\div$, and thus plays no role in the present paper. 

\

\

{\bf 2. Hyperbolic lattice VOA.}

\

Representations of toroidal Lie algebras are best described using the language
of vertex algebras.
Here we sketch the construction of a hyperbolic lattice VOA. The details can be found
in [BBS] or [B3]. We refer to [K] for the definition and properties of vertex algebras.

Consider a hyperbolic lattice $\Hyp$, which is a free abelian group on $2N$
generators 
\break
$\{ u_i , v_i |  i = 1, \ldots, N \}$ with the symmetric bilinear
form
$$ ( \cdot | \cdot ) : \quad \Hyp \times \Hyp \rightarrow \Z ,$$
defined by 
$$ (u_i | v_j) = \delta_{ij} , \quad (u_i | u_j) = (v_i | v_j) = 0.$$
Note that the form $(\cdot | \cdot)$ is non-degenerate and $\Hyp$ is an 
even lattice, i.e., $(x | x) \in 2\Z$.

Consider an infinite-dimensional Heisenberg algebra 
$$\widehat H = \left( \C[t, t^{-1}] \ot_{\Z} \Hyp \right) \oplus \C K$$
with the bracket
$$ [x(n), y(m)] = n (x | y) \delta_{n, -m} K, \quad x,y \in \Hyp, 
\quad [\widehat H, K] = 0. $$
Here and in what follows, we are using the notation $x(n) = t^n \ot x$.
The natural $\Z$-grading of $\widehat H$ yields the decomposition
${\widehat H} = {\widehat H}_- \oplus {\widehat H}_0 \oplus {\widehat H}_+$.

 The hyperbolic lattice VOA $V_\Hyp$ is a tensor product of a twisted group
algebra $\C[\Hyp]$ of the lattice and the Fock space $S(\widehat H_-)$.

The elements of $\widehat H$ of degree zero act on $V_\Hyp$ as follows:
$$ x(0) e^y = (x | y) e^y, \quad K e^y = e^y. $$

The Virasoro element in $V_\hyp$ is $\omega_\hyp = \sum\limits_{p=1}^N
u_p (-1) v_p(-1) \o$, where $\o = e^0$ is the identity element
of $V_\hyp$. The rank of $V_\hyp$ is $2N$.

 The $Y$-map is defined on the basis elements of $\C[\Hyp]$ by
$$ Y(e^x, z) = \exp \left( \sum\limits_{j \geq 1} {x(-j) \over j} z^j \right)
  \exp \left( - \sum\limits_{j \geq 1} {x(j) \over j} z^{-j} \right)
e^x z^x ,$$
where $z^x e^y = z^{(x|y)} e^y$.

For the generators of the Fock space, the $Y$-map is given by
$$ Y(u_p (-1) \o, z) = u_p(z) = \sum\limits_{j\in \Z} u_p(j) z^{-j-1}, \quad
Y(v_p (-1) \o, z) = v_p(z) = \sum\limits_{j\in \Z} v_p(j) z^{-j-1} .$$

In the construction of the toroidal VOAs we would need not $V_\hyp$ itself,
but its sub-VOA $V_\hyp^+$:
$$  V_\hyp^+ = S(\widehat H_-) \ot \C[\Hyp^+] ,$$
where $\Hyp^+$ is the isotropic sublattice of $\Hyp$ 
generated by $\{ u_i |  i = 1, \ldots, N \}$. Here $\C[\Hyp^+]$ is the ordinary
group algebra of $\Hyp^+$. 

The Virasoro element of $V_\hyp^+$ is the same as in $V_\hyp$, and so
the rank of $V_\hyp^+$ is also $2N$.

We will need a class of modules for $V_\hyp^+$.
Fix $\alpha \in \C^N, \beta\in\Z^N$. Then the space
$$ M_\hyp^+(\alpha,\beta) = S(\widehat H_-) \ot 
e^{\alpha u + \beta v } \C[\Hyp^+] $$ has a structure of 
an irreducible
VOA module for $V_\hyp^+$ [BBS], [B3]. Here we are using the notations
$\alpha u = \alpha_1 u_1 + \ldots +  \alpha_N u_N$, etc.

\

\

{\bf 3. VOA associated with the twisted Virasoro-affine algebra.}

\

It has been shown in [B3] that the representation theory of the full
toroidal Lie algebra is controlled by a VOA which is 
a tensor product of two VOAs -- $\VH$ described above 
and a twisted Virasoro-affine VOA $V_\sf (\gamma)$
which we are going to discuss now.

 Let $\df$ be a finite-dimensional reductive Lie algebra. 
 
Consider the 
semi-direct product of the Lie algebra of vector fields on a circle 
with a loop algebra:
$$ \tf = \Der \C[t_0, t_0^{-1}] \ltimes  \left( \C[t_0, t_0^{-1}] \otimes \df \right) .$$
The Lie algebra $\tf$ is a perfect Lie algebra, and we consider its universal central extension
$\f$, which we call the {\it twisted Virasoro-affine algebra} corresponding to $\df$.

In case when $\df$ is a trivial one dimensional Lie algebra, the structure of the universal central
extension of $\tf$ was determined in [ACKP], Proposition 2.1.

 In the present paper, we are interested in the case when $\df = \dg \oplus \glN$, which we fix for
the rest of the paper. Using the well-known description of the universal central extensions for
the loop algebras together with Theorem 1.5.2 of [F], one can modify the argument of Proposition
2.1 in [ACKP] to show that for $\df = \dg \oplus \glN$, the corresponding twisted Virasoro-affine
algebra $\f$ is a 5-dimensional extension of $\tf$:
$$ \f = \tf \oplus \C C_\sdg \oplus \C C_\sln \oplus \C C_\Hei \oplus \C C_\HV \oplus \C C_\Vir .$$ 
Let us describe this central extension in more detail.

The Lie algebra $\f$ contains four subalgebras -- a Virasoro algebra 
$\Vir = \Der \C[t_0, t_0^{-1}] \oplus  \C C_\Vir$, with the bracket
$$[L(n), L(m)] = (n-m) L(n+m) + {n^3 - n \over 12} \delta_{n,-m}  C_\Vir,$$
two affine algebras,
$\wdg = \C [t_0, t_0^{-1}] \otimes \dg \oplus \C C_\sdg$ and 
$\wsl = \C [t_0, t_0^{-1}] \otimes \slN \oplus \C C_\sln$,
and an infinite-dimensional Heisenberg algebra 
$\Hei = \C [t_0, t_0^{-1}] \otimes I \oplus \C C_\Hei$, where $I$ is the identity matrix in $\glN$,
and the bracket is given by
$$[I(n), I(m)] = n \delta_{n,-m} C_\Hei .$$
The action of the Virasoro algebra on affine subalgebras is the usual one, whereas
for the Heisenberg algebra it is twisted with a cocycle:
$$ [L(n), I(m)] = -m I(m+n) - (n^2 + n) \delta_{n,-m}  C_\HV. $$ 
Here $L(n) = -t_0^{n+1} {d \over dt_0}$, $I(m) = t_0^m \otimes I$.

The Virasoro and the Heisenberg subalgebras together generate a subalgebra $\HVir$ in $\f$
called the twisted Virasoro-Heisenberg algebra:
$$ \HVir = \Der \C[t_0, t_0^{-1}] \oplus \left( \C[t_0, t_0^{-1}] \otimes I \right)
 \oplus \C C_\Hei \oplus \C C_\HV \oplus \C C_\Vir .$$
 
 In [B3] we have constructed a VOA $V_\sf (\gamma)$ associated with the twisted Virasoro-affine
algebra $\f$ and a central character $\gamma$: 
$$\gamma(C_\sdg) = c_\sdg, \quad
\gamma(C_\sln) = c_\sln, \quad
\gamma(C_\Hei) = c_\Hei, \quad
\gamma(C_\HV) = c_\HV   , \quad
\gamma(C_\Vir) = c_\Vir .$$ 
For a finite-dimensional irreducible $\dg$-module $V$, a finite-dimensional  
irreducible $\slN$-module $W$, and two constants $h_\Vir, h_\Hei \in \C$,
we get a generalized Verma module $M_\sf (V,W,h_\Hei,h_\Vir,\gamma)$ and its
irreducible quotient $L_\sf (V,W,h_\Hei,h_\Vir,\gamma)$ as modules for both $\f$ and 
the VOA $V_\sf (\gamma)$ (for details see section 3.4 in [B3]).

Let $\wVir$ be another copy of the Virasoro algebra with the basis 
$\{ \wL (n), \wCL \, | \, n\in\Z \}$.
Consider a  semidirect product $\bar \f$ of $\wVir$
with affine algebras $\wdg \oplus \wsl$.
We want to distinguish between $\Vir$ and $\wVir$ because we will be using non-trivial
embeddings of $\bar \f$ into $\f$ (Lemma \emb \ below).
A similar construction applied to $\bar \f$ instead of $\f$ results in a VOA
$V_\sbarf (c_\sdg, c_\sln, {\bar c}_\Vir)$,
\break
a generalized Verma module $M_\sbarf (V,W,{\bar h}_\Vir, c_\sdg, c_\sln, {\bar c}_\Vir)$
and the irreducible module
\break
$L_\sbarf (V,W,{\bar h}_\Vir, c_\sdg, c_\sln, {\bar c}_\Vir)$ for this VOA. 

We will also be using
the affine vertex algebras $V_\swdg (c_\sdg)$, $V_\wsl (c_\sln)$, the
Virasoro VOA $V_\Vir (c_\Vir)$, twisted Virasoro-Heisenberg VOA 
$V_\HVir (c_\Hei, c_\HV, c_\Vir)$,
their Verma modules 
\break
and their irreducible highest weight modules
$L_\swdg (V, c_\sdg)$, $L_\wsl (W, c_\sln)$,
$L_\Vir (h_\Vir, c_\Vir)$ and 
\break
$L_\HVir (h_\Hei, h_\Vir, c_\Hei, c_\HV, c_\Vir)$.

Applying the well-known Sugawara construction, we can decompose the VOA
\break
$V_\sbarf (c_\sdg, c_\sln, {\bar c}_\Vir)$ and its irreducible modules into
tensor products:

{\bf Proposition \FA.}
{\it
 Let $c_\sdg \neq - h^{\vee}, c_\sln \neq -N$, where $h^{\vee}$ is 
the dual Coxeter number for $\wdg$.
Then 

(i) the VOA $V_\sbarf (c_\sdg, c_\sln, {\bar c}_\Vir)$ decomposes into a tensor product of three VOAs:
$$V_\sbarf (c_\sdg, c_\sln, {\bar c}_\Vir) \cong V_\swdg (c_\sdg) \otimes V_\wsl (c_\sln)  
\otimes V_\Vir (c^\prime_\Vir) ,$$
where 
$$c^\prime_\Vir = {\bar c}_\Vir - {c_\sdg \dim (\dg) \over c_\sdg
+ h^\vee} - {c_\sln (N^2 -1) \over c_\sln + N}, \eqno{(\cvp)}$$

(ii) the irreducible highest weight $\barf$-module
$L_\sbarf (V,W,{\bar h}_\Vir, c_\sdg, c_\sln, {\bar c}_\Vir)$
decomposes into a tensor product of irreducible highest weight
modules for the affine algebras $\wdg$, $\wsl$
and the Virasoro modules:
$$L_\sbarf (V,W,{\bar h}_\Vir, c_\sdg, c_\sln, {\bar c}_\Vir)
\cong 
L_{\swdg} (V, c_\sdg) \otimes L_{\wsl} (W, c_\sln) 
\otimes L_{\Vir} (h_\Vir^\prime, c_\Vir^\prime),$$
where $c_\Vir^\prime$ is given by (\cvp) and
$$h_\Vir^\prime = {\bar h}_\Vir - {\Omega_V \over 2(c_\sdg + h^\vee)} - {\Omega_W \over 2(c_\sln + N)}. \eqno{(\hvp)}$$

Here $\Omega_V$ and $\Omega_W$ are the eigenvalues of the Casimir operators of $\dg$ and $\slN$ on $V$
and $W$ respectively. 
}

\

\

{\bf 4. Irreducible modules for the full toroidal Lie algebra.}

\

 In [B3] we introduced a category of bounded modules for the full toroidal Lie algebras
and described irreducible modules in that category. When we restrict these modules
to the subalgebra $\gdiv$, they become reducible. 
The goal of this paper is to describe a class of irreducible modules for $\gdiv$ that
occur as quotients in this reduction.

 We begin by recalling a result of [B3]:
 
{\bf Theorem 4.1. [B3]}. 
{\it
Let $c \neq 0$. 
Let $\VH$ be a sub-VOA of the hyperbolic lattice VOA 
and let $V_\sf (\gamma_0)$ be the enveloping vertex algebra for the twisted Virasoro-affine
Lie algebra $\f$ with $\df = \dg \oplus \glN$, where the central character $\gamma_0$ given
by the following values: 
$$ c_\sdg = c, \quad \quad c_\sln = 1 - \mu c,  \quad\quad c_\Hei = N(1-\mu c),$$
$$  c_{\HV} = {N \over 2},  \quad\quad c_\Vir = 12  \mu c  - 2N. \eqno{(\gamo)}$$ 
Then the VOA 
$$\VH \otimes V_\sf (\gamma_0)$$
has a structure of a module over the full toroidal Lie algebra $\g(\mu)$
with the action given by 
$$ k_0 (r,z) =
\sum_{j\in\Z} t_0^j t^r k_0 z^{-j} \mapsto c Y(\eru, z) , \quad r\in\Z^N,  \eqno{(\ko)}$$
$$ k_a (r,z) =
\sum_{j\in\Z} t_0^j t^r k_a z^{-j-1} \mapsto c u_a(z) Y(\eru, z), \eqno{(\ka)}$$ 
$$ g (r,z) =
\sum_{j\in\Z} t_0^j t^r g z^{-j-1} \mapsto g(z) Y(\eru, z), \quad g \in \dg, \eqno{(\gr)}$$ 
$$ d_a (r,z) =
\sum_{j\in\Z} t_0^j t^r d_a z^{-j-1} \mapsto  :v_a(z) Y(\eru, z): 
+ \sum_{p=1}^N r_p E_{pa}(z) Y(\eru,z) , \eqno{(\da)}$$ 
$$ \td_0 (r,z) = 
\sum_{j\in\Z} t_0^j t^r \td_0 z^{-j-2} \mapsto \hbox{\hskip 9cm}$$
$$  :\omega(z) Y(\eru, z): 
+ \sum_{p,s=1}^N r_p u_s(z) E_{ps}(z) Y(\eru,z) 
+ (\mu c -1) \sum_{p=1}^N r_p \left( {\d \over \d z} u_p(z) \right) Y(\eru,z)  . \eqno{(\tdd)}$$
Here 
$$ t_0^j t^r \td_0 = - t_0^j t^r d_0 + \mu (j+{1\over 2}) t_0^j t^r k_0 $$
and $\om (z)$ is the Virasoro field in the VOA $\VH \otimes V_\sf (\gamma)$, which has
rank $12 \mu c$. In (\gr)-(\tdd) above we are using the notations
$g(z) = Y(g(-1) \o, z)$, $E_{ab} (z) = Y(E_{ab} (-1) \o, z)$, where $E_{ab}$ are the basis
elements of $\glN$.
}
 
 Taking simple modules for the VOAs in the previous theorem, we obtain 
irreducible representations of $\g (\mu)$:

{\bf Theorem 4.2. [B3]}. 
{\it
Let $M_\Hyp^+ (\alpha, \beta)$ be a simple module for
$\VH$ and let 
\break
$M_\sf$ be a VOA module for $V_\sf (\gamma_0)$. Then
$$ M_\Hyp^+ (\alpha, \beta) \otimes M_\sf$$
is a module for the full toroidal Lie algebra $\g(\mu)$ with
the action given by the formulas (\ko)-(\tdd).
This module is irreducible as a $\g(\mu)$-module whenever $M_\sf$ is an irreducible
$\f$-module.
}

\

\

{\bf 5. Irreducible representations for $\gdiv$.}
 
\

Now we will study the restriction of the modules for the toroidal Lie algebra to
the subalgebra $\gdiv (\mu) = (\R \otimes \dg) \oplus \K \oplus \Ddiv$. This subalgebra is
spanned by the elements $t_0^j t^r k_0$, $t_0^j t^r k_p$, $t_0^j t^r g$,
$d_0$, $t_0^j d_p$,
$$ r_b t_0^j t^r d_a - r_a t_0^j t^r d_b, \eqno{(\drab)} $$
and
$$ t_0^j t^r \wda = r_a t_0^j t^r \td_0 + j t_0^j t^r d_a 
+ {r_a \over 2cN} (N-1+\mu c) t_0^j t^r k_0 . \eqno{(\wdoa)}$$
The reason for adding the last term in (\wdoa) will become clear later.

The elements $t_0^j t^r k_0$, $t_0^j t^r k_p$, $t_0^j t^r g$
correspond to the fields (\ko)-(\gr), while
$$d_p(z) = \sum\limits_{j\in\Z} t_0^j d_p z^{-j-1} \mapsto v_p (z). \eqno{(\dpz)} $$
Collect the elements of the form (\drab) and (\wdoa) into the fields:
$$d_{ab} (r,z) = \sum_{j\in\Z} 
\left( r_b t_0^j t^r d_a - r_a t_0^j t^r d_b \right) z^{-j-1},  $$
$$\wda (r, z) = \sum_{j\in\Z} t_0^j t^r \wda z^{-j-2}.$$
It follows from (\da) that $d_{ab} (r,z)$ acts on $\VH \otimes V_\sf (\gamma_0)$ by
$$ d_{ab} (r,z) \mapsto
: \left( r_b v_a(z) - r_a v_b(z) \right) Y(\eru,z) :$$
$$ + r_b \sum\limits_{{p=1} \atop {p\neq a}}^N r_p E_{pa} (z) Y(\eru,z)
- r_a \sum\limits_{{p=1} \atop {p\neq b}}^N r_p E_{pb} (z) Y(\eru,z)
+ r_a r_b \left( E_{aa} - E_{bb} \right) (z) Y(\eru,z). \eqno{(\Ydab)} $$
 It is easy to see that the moments of the fields 
(\ko)-(\gr), (\dpz), (\Ydab) generate $V_\hyp^+, V_\swdg$ and $V_\wsl$
(see the proof of Theorem 5.1 in [BBS]). 

Using (\ko)-(\tdd), we obtain that  
$\wda (r, z)$ is represented on $\VH \otimes V_\sf (\gamma_0)$
in the following way:
$$ \wda (r, z) \mapsto r_a Y \left( \om_{(-1)} \eru + 
\sum_{p,s=1}^N r_p u_s(-1) \eru \ot E_{ps} (-1) 
+ \left(\mu c - 1 \right) \sum_{p=1}^N r_p u_p(-2) \eru, z \right)$$
$$- \left( z^{-1} + {\d \over \d z} \right) Y \left( v_a(-1) \eru + 
\sum_{p=1}^N r_p \eru \ot E_{pa} (-1) , z \right) 
+ {r_a z^{-2} \over 2N} (N-1+\mu c) Y \left( \eru, z \right). \eqno{(\dwd)}$$

Since $V_\hyp^+ \ot V_\swdg \ot  V_\wsl$ is generated by $\gdiv$, we shall consider
in (\dwd) only the terms that involve $V_\HVir$, together with the last summand:
$$r_a Y \left( \eru \ot \om_\HVir, z \right) + {r_a \over N}
Y \left( \sum_{p=1}^N r_p u_p(-1) \eru \ot I(-1), z \right)$$
$$- {r_a \over N} Y \left( D(\eru \ot I(-1)), z \right)
- z^{-1} {r_a \over N} Y \left( \eru \ot I(-1), z \right)
+ {r_a z^{-2} \over 2N} (N-1+\mu c) Y \left( \eru, z \right)$$
$$ = r_a \left( Y \left( \eru \ot \om_\HVir, z \right)
- {1 \over N} Y \left( \eru \ot D(I(-1)), z \right) \right.$$
$$\left.
- z^{-1} {1 \over N} Y \left( \eru \ot I(-1), z \right)
+ {z^{-2} \over 2N} (N-1+\mu c) Y \left( \eru, z \right) \right)$$
$$ = r_a Y \left(\eru, z \right) \ot
\left( Y \left(\om_\HVir, z \right)
- {1 \over N} \left(z^{-1} + {\d \over \d z} \right) Y \left( I(-1), z \right)
+ {z^{-2} \over 2N} (N-1+\mu c) \Id \right) .$$
In the above calculation we used the fact that 
$D(\eru) = \mathop\sum\limits_{p=1}^N r_p u_p (-1) \eru$ in $\VH$.

To understand the structure of the expression 
$$Y \left(\om_\HVir, z \right)
- {1 \over N} \left(z^{-1} + {\d \over \d z} \right) 
Y \left( I(-1), z \right)
+ {z^{-2} \over 2N} (N-1+\mu c) \Id, \eqno{(\expr)}$$
we consider the following

{\bf Lemma \emb.} 
{\it
Let $\wVir$ be the Virasoro algebra with the basis
$\left\{ \wL (n), \wCL \right\}$. For any $\sigma \in\C$ the map
$$\rho_\sigma : \quad \wVir \rightarrow \HVir,$$
given by 
$$\rho_\sigma (\wL (n)) = L(n) + \sigma n I(n) 
+ \delta_{n,0} (\sigma C_\HV - {\sigma^2 \over 2} C_\Hei),$$
$$\rho_\sigma (\wCL) = C_\Vir + 24 \sigma C_\HV -12 \sigma^2 C_\Hei,$$
is an embedding of Lie algebras.
}

\

{\bf Corollary \Memb.}
{\it
(i) Let $h_\Vir, h_\Hei, c_\Vir, c_\HV, c_\Hei \in \C$, and let
$$\bar h_\Vir = h_\Vir + \sigma c_\HV - {\sigma^2 \over 2} c_\Hei, \quad
\bar c_\Vir = c_\Vir + 24 \sigma c_\HV - 12 \sigma^2 c_\Hei . \eqno{(\barhc)}$$
 The homomorphism
$\rho_\sigma$ extends to the embedding of the Verma module
$M_\wVir (\bar h_\Vir, \bar c_\Vir)$ for the Virasoro algebra $\wVir$, 
into the Verma module 
$M_\HVir (h_\Hei, h_\Vir, c_\Hei, c_\HV, c_\Vir)$ for the twisted Virasoro-Heisenberg algebra $\HVir$.

\noindent
(ii) The homomorphism $\rho_\sigma$ also extends to the embedding of the 
Lie algebras 
$$\rho_\sigma : \barf \rightarrow \f ,$$
with the identical action on the subalgebras $\wdg$ and $\wsl$. This yields an embedding of the generalized
Verma modules
$$M_\sbarf (V,W,{\bar h}_\Vir, c_\sdg, c_\sln, {\bar c}_\Vir) \arsubset 
M_\sf (V,W,h_\Hei,h_\Vir,\gamma) .$$

\noindent
(iii) Let $c_\HV = {N \over 2}, c_\Hei = N (1 - \mu c), \sigma = {1 \over N}$. 
Under the map $\rho_\sigma$
we have the correspondence of the fields 
$$\sum_{n\in\Z} \wL(n) z^{-n-2} \mapsto \sum_{n\in\Z} L(n) z^{-n-2}
- {1 \over N} \left( z^{-1} + {\d \over \d z} \right) \sum_{n\in\Z} I(n) z^{-n-1}
+ z^{-2}\left( {1 \over 2} - {1 - \mu c \over 2N} \right) \Id . \eqno{(\embd)}$$
}

The proof of the Lemma is a straightforward computation, and the Corollary is an
immediate consequence. Note that the field in the right hand side of (\embd)
coincides precisely with (\expr). Thus when $c_\HV = {N \over 2}$,
$c_\Hei = N (1 - \mu c)$,
(as we have in Theorem 4.1), the components of this field satisfy the relations
of the Virasoro algebra with the value of the central charge $\wcL = c_\Vir + 12
- {12 \over N} (1 - \mu c)$.

{\bf Remark.} Note that the map $\rho_\sigma$ does not extend to the homomorphism
of VOAs, since $\wL(-1) \o = 0$, while $\rho_\sigma (\wL (-1)) \o \neq 0$ for $\sigma \neq 0$.

 In order to construct representations, one needs the principle of 
``preservation of identities'' for the vertex algebras ([Li], Lemma 2.3.5). In our
situation, however,
the field (\expr) involves vertex operators that are shifted by powers of $z$.
To deal with such expressions we need to establish a generalization of this principle
(Lemma \presr \ below).

 A convenient set-up for working with the fields shifted by powers of $z$ is given 
by the construction of the affinization of a vertex algebra. Let us review this construction.

 We recall that for a commutative associative unital algebra $R$ with a derivation $D$,
one can define a vertex algebra structure on $R$ by the formula
$$ Y(a,z) b = \sum\limits_{n=0}^\infty {z^n \over n!} D^n (a) b .$$
Let us now take $R$ to be the algebra of Laurent polynomials $\C [t,t^{-1}]$ with the
derivation $D = {d \over dt}$. This vertex algebra has a 1-dimensional module $Z$
such that 
$$ Y_Z (t^k, z) = z^k \Id_Z .$$
The affinization of a vertex algebra $V$ is defined as a tensor product of vertex algebras
$V \otimes \C[t,t^{-1}]$. Let $M$ be a module for the vertex algebra $V$. Then we can view
$M \cong M \otimes Z$ as a module for the affinization of $V$:
$$Y_{M\otimes Z} (a \otimes t^k, z) = z^k Y_M (a,z) .$$

{\bf Lemma \presr.}
{\it
Let $V$ be a VOA and let $M$ be a VOA module for $V$. 
Let $a \in V \otimes \C[t,t^{-1}]$.

(i) If $M$ is a faithful VOA module for $V$ and $Y_{M\otimes Z} (a,z) = 0$ then $a=0$.

(ii) For $a, b, c^0, \ldots c^N \in V \otimes \C[t,t^{-1}]$, if
$$ \left[ Y_{V\otimes Z} (a,z), Y_{V\otimes Z} (b,w) \right] =
\sum\limits_{n=0}^N Y_{V\otimes Z} (c^n,w) \left[ z^{-1} \left( \d \over \d w \right)^n 
\delta \left( w \over z \right) \right] , \eqno{(\relV)}$$
then
$$ \left[ Y_{M\otimes Z} (a,z), Y_{M\otimes Z} (b,w) \right] =
\sum\limits_{n=0}^N Y_{M\otimes Z} (c^n,w) \left[ z^{-1} \left( \d \over \d w \right)^n 
\delta \left( w \over z \right) \right] . \eqno{(\relM)}$$

(iii) If $M$ is a faithful VOA module for $V$ then (\relM) implies (\relV).
}

{\it Proof.} We begin by proving (i). Let $D_M$ be the infinitesimal translation operator on $M$.
Write $a = \sum\limits_s a^s \otimes t^{-s}$, with $a^s \in V$.
If $\sum\limits_s Y_M (a^s, z) z^{-s} = 0$ then 
$$ 0 = z \sum\limits_s  \left[ D_M, Y_M (a^s, z)  \right] z^{-s} 
= z \sum_s \left( {\d \over \d z} Y_M (a^s, z) \right) z^{-s} $$
$$ = z \sum_s \left( {\d \over \d z} Y_M (a^s, z) \right) z^{-s}
- z {\d \over \d z} \left( \sum_s  Y_M (a^s, z) z^{-s} \right) $$
$$ = - \sum\limits_s Y_M (a^s, z) z {\d \over \d z} (z^{-s}) 
= \sum\limits_s s Y_M (a^s, z) z^{-s} .$$
Repeating this argument, we get that for any $m = 0,1,2, \ldots $
$$ \sum\limits_s s^m Y_M (a^s, z) z^{-s}  = 0 .$$
Since the sum in $s$ is finite, we can apply the Vandermonde determinant argument and
derive that $Y_M (a^s, z) = 0$ for all $s$. By the definition of a faithful module,
this implies that all $a^s = 0$.   

 To prove (ii), we use the commutator formula for the vertex algebras and the basic properties of
the delta-function:
$$\left[ \sum_s z^{-s} Y_V (a^s, z), \sum_k w^{-k} Y_V (b^k,w) \right]$$
$$= \sum_{n,i\geq 0} \sum_{s,k} {1 \over (n+i)!} \pmatrix{ -s \cr i \cr}
 w^{-k-s-i} Y_V (a^s_{(n+i)} b^k, w) 
\left[ z^{-1} \left( {\d \over \d w} \right)^n \delta \left( {w\over z} \right) \right]
\quad \hbox{\rm (all sums finite)}. \eqno{(\relVY)}$$
By Corollary 2.2 from [K], we obtain that for all $n \geq 0$,
$$ Y_{V\otimes Z} (c^n, w) = \sum_j w^{-j} Y_V (c^{n,j},w) = 
\sum_{i\geq 0} \sum_{s,k} {1 \over (n+i)!} \pmatrix{ -s \cr i \cr}
 w^{-k-s-i} Y_V (a^s_{(n+i)} b^k, w).$$   
Since $V$ is a faithful VOA module over itself, we get using part (i) of the Lemma
that
$$ c^{n,j} = \sum_{s,k} \sum_{{i\geq 0}\atop {s+k+i = j}} 
{1 \over (n+i)!} \pmatrix{ -s \cr i \cr} a^s_{(n+i)} b^k . \eqno{(\cab)}$$
However the relation (\relVY) holds in every VOA module M. Taking (\cab) into account,
we see that (\relM) also holds.

The proof for part (iii) is similar. We first see that
$$\sum_{n\geq 0} \sum_j w^{-j} Y_M (c^{n,j},w) 
\left[ z^{-1} \left( {\d \over \d w} \right)^n \delta \left( {w\over z} \right) \right] =$$
$$\sum_{n,i\geq 0} \sum_{s,k} {1 \over (n+i)!} \pmatrix{ -s \cr i \cr}
 w^{-k-s-i} Y_M (a^s_{(n+i)} b^k, w) 
\left[ z^{-1} \left( {\d \over \d w} \right)^n \delta \left( {w\over z} \right) \right] .$$
Again using Corollary 2.2 from [K] and part (i) of the Lemma, we obtain that the 
relation
(\cab) holds in $V$. Thus (\relV) also holds. This completes the proof of the Lemma.

\

Now we have done all the preparatory work and now ready to describe the representations 
for $\gdiv$.
  
{\bf Theorem \maind.}
{\it Let $c\neq 0$. Let $\barf$ be the semidirect product of the Virasoro algebra with
the affine Lie algebra $\wdg \oplus \wsl$, and let
$V_\sbarf (c_\sdg, c_\sln, \bar c_\Vir)$ be the enveloping VOA for $\barf$ with the
central charges
$$c_\sdg = c, \quad c_\sln = 1 - \mu c, \quad
\bar c_\Vir = 12 \left( 1 - {1 \over N} \right) + 12 \mu c \left( 1 + {1 \over N} \right) - 2N . \eqno{(\cenc)}$$
Then 

(i) the VOA 
$V_\sgdiv = \VH \otimes V_\sbarf (c_\sdg, c_\sln, \bar c_\Vir)$
is a module for the Lie algebra $\gdiv (\mu)$. 
The action of the fields
$k_0 (r,z), k_a(r,z), g(r,z)$ is given by
the formulas (\ko)-(\gr).
The action of $d_p(z)$ is given by (\dpz).
The field $d_{ab} (r,z)$ acts according to (\Ydab).
The action of $d_0$ is given by 
$$ d_0 \mapsto  - \om_{(1)}, \eqno{(\doo)}$$
where $\om$ is the Virasoro element of the VOA $V_\sgdiv$.
Finally, the field $\wda(r,z)$ is represented by 
$$ r_a Y \left(\om_{(-1)} \eru + \sum_{p,s =1}^N r_p u_s (-1) \eru \ot 
\psi (E_{ps}) (-1) 
+ \left(\mu c - 1 \right) \sum_{p=1}^N r_p u_p(-2) \eru, z \right) $$
$$ - \left( z^{-1} + {\d \over \d z} \right)
Y \left( v_a(-1) \eru + \sum_{p=1}^N r_p \eru \ot 
\psi (E_{pa}) (-1), z \right),  \eqno{(\wdaa)}$$
where $\psi$ is the natural projection $\gln \rightarrow \sln$,
\quad $\psi(X) = X - {1\over N} \tr (X) I$.

(ii) Let $M_\hyp^+(\alpha, \beta)$ be a VOA module for the sub-VOA $V_\hyp^+$
of the hyperbolic lattice VOA, and let $M_\sbarf$ be a module 
for the VOA $V_\sbarf (c_\sdg, c_\sln, \bar c_\Vir)$.
Then 
$$M_\hyp^+(\alpha, \beta) \otimes M_\sbarf$$
is a module for the Lie algebra $\gdiv (\mu)$ with the action transferred 
from $V_\sgdiv$.
}

{\bf Remark.} Since degree zero derivations $d_0$, $d_p$ do not belong to the commutant of $\g_\div$,
their action may be modified by adding arbitrary multiples of the identity operator.
  
{\it Proof.} The Lie bracket in $\gdiv$ may be encoded in the commutator
relations between the fields $k_0 (r,z)$, $k_a(r,z)$, $g(r,z)$, 
$d_p(z)$, $d_{ab} (r,z)$, $\wda(r,z)$,
and the element $d_0$.
We need to show that the same commutator relations hold for their images
(\ko)-(\gr), (\dpz), (\Ydab), (\doo) and (\wdaa). 
It is easy to see that the relations involving $d_0$ in the left hand sides 
hold since $\omega_{(1)}$ acts on the VOA as a degree derivation.
Also, $d_0$ does not belong to the commutant of $\gdiv$ and will not appear
in the right hand sides of the commutator relations.

Actually, there is no need to write down the commutator relations involving the remaining fields
explicitly. All we need to know is that these relations are of the form (\relV),
which follows from the corresponding result for the full toroidal Lie algebra
(Theorem 4.1 in [B3]).
Our strategy is to embed one of the modules for $V_\sgdiv$ into a module for
the full toroidal Lie algebra $\g(\mu)$. This embedding will have the property
that the restriction of the action of $\g(\mu)$ to subalgebra $\gdiv(\mu)$ will coincide
with (\ko)-(\gr), (\dpz), (\Ydab) and (\wdaa).
This will imply that the necessary commutator relations hold in the chosen
module for $V_\sgdiv$. Since the module that we will consider will be faithful,
then by preservation of identities (Lemma \presr), the same required relations
will hold in all VOA modules for $V_\sgdiv$.

Let us carry out this plan. 
Let $h_\Vir = -{1 \over N} c_\HV + {1\over 2N^2} c_\Hei$,
$h_\Hei = 0$, and $c_\sdg$, $c_\Hei$, $c_\HV$,
$c_\sln$, $c_\Vir$, $\wcL$ given by (\gamo), (\barhc)
and let $V$,$W$ be trivial 1-dimensional modules for $\dg$ and $\slN$.
Consider the embedding  given by Corollary \Memb(ii) \ with
$\sigma = {1 \over N}$ 
of the generalized Verma $\bar \f$-module into the generalized Verma module for $\f$:
$$M_\sbarf(V,W,0,c_\sdg, c_\sln, {\bar c}_\Vir) \arsubset 
M_\sf (V,W,0,h_\Vir,\gamma_0) .$$
Under this homomorphism
$$ \wL (z) = \sum_{n\in\Z} \wL (n) z^{-n-2} \mapsto
L(z) - {1 \over N} \left( z^{-1} + {\d \over \d z} \right) I(z) + 
z^{-2} \left( {1\over 2} - {1-\mu c \over 2N} \right) \Id.$$
This map extends to the embedding
$$ V_\hyp^+ \ot M_\sbarf(V,W,0,c_\sdg, c_\sln, {\bar c}_\Vir) \arsubset 
V_\hyp^+ \ot M_\sf (V,W,0,h_\Vir,\gamma_0) .$$
By Theorem 4.2, the latter is a module for the full toroidal Lie algebra 
$\g(\mu)$. We consider the restriction of this representation
to the subalgebra $\gdiv(\mu)$ and claim that
$ V_\hyp^+ \ot M_\sbarf(V,W,0,c_\sdg, c_\sln, {\bar c}_\Vir)$ 
is invariant under the action of $\gdiv(\mu)$. The action of 
$k_0 (r,z)$, $k_a(r,z)$, $g(r,z)$, $d_p(z)$ and $d_{ab}(r,z)$ is given by
(\ko)-(\gr), (\dpz), (\Ydab) and the invariance with respect to these fields is clear. 
Let us show that the action of $\wda(r,z)$ on 
$ V_\hyp^+ \ot M_\sbarf(V,W,0,c_\sdg, c_\sln, {\bar c}_\Vir)$
coincides with (\wdaa).
Indeed, following the computations (\dwd)--(\embd), we get:
$$\wda(r,z) \mapsto 
r_a Y(\eru, z) \otimes \left(
Y \left(\om_\sf, z \right)
- {1 \over N} \left(z^{-1} + {\d \over \d z} \right) 
Y \left( I(-1), z \right)
+ {z^{-2} \over 2N} (N-1+\mu c) \Id \right)$$
$$ + r_a Y \left( \om_\Hyp {}_{(-1)} \eru
+ \sum_{p,s = 1}^N r_p u_s(-1) \eru \ot \psi(E_{ps}) (-1)
+ \left(\mu c - 1 \right) \sum_{p=1}^N r_p u_p(-2) \eru, z \right) $$
$$ - \left( z^{-1} + {\d \over \d z} \right) 
Y \left( v_a (-1) \eru + \sum_{p=1}^N r_p \eru \ot \psi (E_{pa}) (-1), z \right)$$
$$ = r_a Y \left(\left( \om_\Hyp + \rho_\sigma(\om_\sbarf) \right)_{(-1)} \eru, z \right)$$
$$+ r_a Y \left( \sum_{p,s = 1}^N r_p u_s(-1) \eru \ot \psi (E_{ps}) (-1)
+ \left(\mu c - 1 \right) \sum_{p=1}^N r_p u_p(-2) \eru, z \right) $$
$$ - \left( z^{-1} + {\d \over \d z} \right) 
Y \left( v_a (-1) \eru + \sum_{p=1}^N r_p \eru \ot \psi (E_{pa})(-1), z \right),$$
which is the same as (\wdaa). 
Thus the specified action defines a representation of $\gdiv(\mu)$ on 
$ V_\hyp^+ \ot M_\sbarf(V,W,0,c_\sdg, c_\sln, {\bar c}_\Vir)$,
and the fields (\ko)-(\gr), (\dpz), (\Ydab) and (\wdaa) satisfy the
relations that reflect the Lie bracket in $\gdiv(\mu)$.
This module is a faithful VOA module for $V_\sgdiv$, since $V_\sgdiv$ itself is its
factor module. Thus by the preservation of identities, Lemma \presr,
the required commutator
relations hold in $V_\sgdiv$ and in all VOA modules for $V_\sgdiv$. This establishes
the claim of the theorem.

\

In the next theorem we give the description of the irreducible modules for $\gdiv$.

\

{\bf Theorem \irrd.} 
{\it
Let $c \neq 0$ and let the constants $c_\sdg, c_\sln, \bar c_\Vir$ be given by (\cenc).

(i) Let $V$ be a finite-dimensional irreducible $\dg$-module, $W$ be 
a finite-dimensional irreducible $\slN$-module, ${\bar h}_\Vir \in \C$ and let
$L_\sbarf (V,W,{\bar h}_\Vir,c_\sdg, c_\sln, \bar c_\Vir)$ be an irreducible highest weight
module for the Lie algebra $\barf$.
Let $\alpha\in\C^N, \beta\in\Z^N$, and let
$M_\hyp^+ (\alpha,\beta)$ be the irreducible VOA module for $V_\hyp^+$.
Then
$$ L_\sgdiv = M_\hyp^+ (\alpha,\beta) \otimes L_\sbarf (V,W,{\bar h}_\Vir,c_\sdg, c_\sln, \bar c_\Vir)$$
is an irreducible module for the Lie algebra $\gdiv(\mu)$ with the action given by (\ko)-(\gr),
(\dpz), (\Ydab), (\doo) and (\wdaa).

(ii) If, in addition, $c_\sdg \neq -h^\vee$, where $h^\vee$ is the dual Coxeter number of
$\wdg$, and $c_\sln = 1 - \mu c \neq -N$ then
$$ L_\sgdiv \cong M_\hyp^+ (\alpha,\beta) \otimes L_\swdg (V, c_\sdg) \otimes L_\wsl (W, c_\sln)
\otimes L_\Vir(h_\Vir^\prime, c_\Vir^\prime) ,$$
where $h_\Vir^\prime$ is given by (\hvp) and
$$c_\Vir^\prime = 12 \left( 1 - {1 \over N}\right) + 12 \mu c  \left( 1 + {1 \over N}\right) - 2N 
- {c \dim(\dg) \over c + h^\vee} - {(1-\mu c)(N^2-1) \over 1 - \mu c + N}. 
\eqno{(\cencp)}$$
}

The proof of this theorem is completely analogous to the proof of 
Theorem 5.3 in [B3] and will be omitted.

\

We conclude the paper with the following observation. 
Let us now set $\dg = (0)$. Then Theorem \irrd \ gives a family of irreducible
representations for the Lie algebras
$$\DK (\mu) =  \Ddiv \oplus \K.$$
These Lie algebras possess a non-degenerate symmetric invariant bilinear form
and may be viewed as extended affine Lie algebras of rank 0.
In fact, all such algebras with $\mu \neq 0$ are isomorphic to each other.
We will thus normalize $\mu$ to be 1.
If we take trivial modules for the Virasoro and for the affine algebra $\wsl$,
then only when $N=12$ we arrive at the irreducible representation of $ \DK (\mu)$ just
on the lattice part $V^+_\hyp$.
We get the following remarkable result:

{\bf Theorem \TDK .} 
{\it
Let $N = 12$ and let $\mu = c = 1$. Then $V^+_\hyp$ has a structure of an irreducible module for
$\DK (1)$. 
The action of the Lie algebra is given by
$$ k_0 (r, z) \mapsto Y(\eru, z), \quad
k_p (r, z) \mapsto u_p (z) Y(\eru, z), $$
$$d_{ab}(r,z) \mapsto :\left( r_b v_a(z) - r_a v_b(z) \right) Y(\eru,z): ,$$
$$ d_0 \mapsto  \Id - \om_{(1)} , \quad d_p (z) \mapsto v_p(z),  $$
$$\wda (r, z) \mapsto r_a :\om(z) Y(\eru, z): - 
\left( z^{-1} + {\d \over \d z} \right) :v_a(z) Y (\eru, z):.$$
}

The character of this module with respect to the diagonalizable operators 
$d_0$, $d_1$, $\ldots$, $d_N$ has nice modular properties -- 
it is a product of $12$ delta-functions with the 
$-24$-th power of the Dedekind $\eta$-function:
$$ \char V^+_\hyp = q_0 \prod_{k=1}^\infty \left( 1 - q_0^{-k} \right)^{-24} 
\times \prod_{p=1}^{12} \sum_{n \in\Z} q_p^n .$$ 
 
\

\vfill\eject

{\bf References:}

\

\item{[AABGP]}  B.N.~Allison,  S.~Azam, S.~Berman,  Y.~Gao, A.~Pianzola,  
{Extended affine Lie algebras and their root systems.} 
Mem.Amer.Math.Soc. {126}, no. 603, 1997.

\item{[ABFP]} B.~Allison, S.~Berman, J.~Faulkner, A.~Pianzola,
{Realizations of graded-simple algebras as loop algebras.}
math.RA/0511723.

\item{[ACKP]} E.~Arbarello, C.~De Concini, V.G.~Kac, C.~Procesi,
{Moduli spaces of curves and representation theory.}
Commun.Math.Phys. {117}, (1988), 1-36.

\item{[BB]} S.~Berman,  Y.~Billig, 
{Irreducible representations for toroidal Lie algebras.} 
J.Algebra {221}, (1999) 188-231.

\item{[BBS]} S.~Berman, Y.~Billig, J.~Szmigielski, 
{Vertex operator algebras and the representation theory of toroidal algebras.}
in ``Recent developments in infinite-dimensional Lie algebras and conformal
field theory'' (Charlottesville, VA, 2000), 
Contemp. Math. {297}, 1-26, Amer. Math. Soc., 2002.

\item{[BGK]} S.~Berman, Y.~Gao, Y.~Krylyuk, 
{Quantum tori and the structure of elliptic quasi-simple Lie algebras.} 
J.Funct.Analysis {135}, (1996) 339-389.

\item{[B1]} Y.~Billig,  
{Energy-momentum tensor for the toroidal Lie algebras.}
math.RT/0201313. 

\item{[B2]} Y.~Billig, 
{Magnetic hydrodynamics with asymmetric stress tensor.}
J.Math.Phys. {46},  (2005) 13pp.

\item{[B3]} Y.~Billig, 
{A category of modules for the full toroidal Lie algebra.}
math.RT/0509368.

\item{[BL]} Y.~Billig, M.~Lau, 
{Irreducible modules for extended affine Lie algebras.}
in preparation.

\item{[EM]} S.~Eswara Rao,  R.V.~Moody,  
{Vertex representations for $n$-toroidal
Lie algebras and a generalization of the Virasoro algebra.}
Commun.Math.Phys. {159}, (1994) 239-264.

\item{[F]} D.B.~Fuks,
{Cohomology of infinite-dimensional Lie algebras.}
Consultants Bureau, N.Y., 1986.

\item{[K]}  V.~Kac,  
{Vertex algebras for beginners.} 
Second edition, University Lecture Series, {10}, A.M.S., 1998.

\item{[Kas]} C.~Kassel, 
{K\"ahler differentials and coverings of complex simple 
Lie algebras extended over a commutative ring,}
J.Pure Applied Algebra {34}, (1984) 265-275.

\item{[L]} T.A.~Larsson, 
{Lowest-energy representations of non-centrally extended diffeomorphism 
algebras.}
Commun.Math.Phys. {201}, (1999) 461-470.

\item{[Li]} H.~Li, 
{Local systems of vertex operators, vertex superalgebras and modules.}
{J.Pure Appl.Algebra} {109}, (1996) 143-195.

\item{[N]} E.~Neher,
{Extended affine Lie algebras.}
{C.R.Math.Acad.Sci.Soc.R.Can.} {26}, (2004) 90-96.

\item{[VD]} S.M.~Vi{\v{s}}ik, F.V.~Dol{\v z}anski{\u{i}}, 
{Analogs of the Euler-Lagrange equations and magnetohydrodynamics connected with
Lie groups.}
Dokl. Akad. Nauk SSSR, {19}, (1978) 149-153.

\end